\documentclass{amsart}

\newtheorem{theorem}{Theorem}

\newtheorem{proposition}[theorem]{Proposition}

\theoremstyle{definition}

\newtheorem{example}[theorem]{Example}

\theoremstyle{remark}
\newtheorem{remark}[theorem]{Remark}

\numberwithin{equation}{section}
\usepackage{amssymb}

\def\BB{\mathbb B}
\def\CC{\mathbb C}
\def\EE{\mathbb E}
\def\NN{\mathbb N}
\def\RR{\mathbb R}
\def\WW{\mathbb W}
\def\too{\longrightarrow}
\def\wdht{\widehat}
\def\wdtl{\widetilde}

\def\sprod#1#2{\langle #1,#2 \rangle}
\def\CF{\mathcal F}
\def\Vol{\operatorname{Vol}}
\DeclareSymbolFont{EulerScriptBold}{U}{eus}{b}{n}
\DeclareSymbolFontAlphabet\eusb{EulerScriptBold}
\def\EM{\eusb M}
\let\epsilon=\varepsilon
\def\eps{\epsilon}
\def\phi{\varphi}
\def\kappa{\varkappa}
\def\Lambda{\varLambda}
\def\Omega{\varOmega}

\begin{document}

\title
{On the upper semicontinuity of the Wu metric}

\author{Marek Jarnicki}
\address{Jagiellonian University, Institute of Mathematics,
Reymonta 4, 30-059 Krak\'ow, Poland}
\email{jarnicki@im.uj.edu.pl}

\author{Peter Pflug}
\address{Carl von Ossietzky Universit\"at Oldenburg, Institut f\"ur Mathematik,
Postfach 2503, D-26111 Oldenburg, Germany}
\email{pflug@mathematik.uni-oldenburg.de}
\thanks{Both authors were supported in part by KBN grant no.~5 P03A 033 21
and by DFG grant no.~227/8-1.}



\begin{abstract}
We discuss continuity and upper semicontinuity of the Wu
pseudometric.
\end{abstract}

\maketitle

The Wu pseudometric has been introduced by H. Wu in \cite{Wu1993}
(and \cite{Wu1987}). Various properties of the Wu metric may be
found for instance in \cite{CheKim1996}, \cite{CheKim1997}, \cite{Kim1998},
\cite{CheKim2003}, \cite{Juc2003}. Nevertheless, it seems that even quite
elementary properties of this metric are not completely understood, e.g.~its
upper semicontinuity.

First, let us formulate the definition of the Wu pseudometric in
an abstract setting. Let $h:\CC^n\too\RR_+$ be a $\CC$--seminorm.
Put:

$I=I(h):=\{X\in\CC^n: h(X)<1\}$ ($I$ is convex),

$V=V(h):=\{X\in\CC^n: h(X)=0\}\subset I$ ($V$ is a vector subspace of $\CC^n$),

$U=U(h):=$ the orthogonal complement of $V$ with respect to the standard
Hermitian scalar product $\sprod{z}{w}:=\sum_{j=1}^nz_j\overline w_j$ in
$\CC^n$,

$I_0:=I\cap U$, $h_0:=h|_U$ ($h_0$ is a norm, $I=I_0+V$).

For any {\it pseudo--Hermitian scalar product}
$s:\CC^n\times\CC^n\too\CC$, let
$$
q_s(X):=\sqrt{s(X,X)},\; X\in\CC^n,\quad \EE(s):=\{X\in\CC^n: q_s(X)<1\}.
$$
Consider the family $\CF$ of all pseudo--Hermitian scalar products
$s:\CC^n\times\CC^n\too\CC$ such that $I\subset\EE(s)$, equivalently,
$q_s\leq h$. In particular,
$$
V\subset I= I_0+V\subset\EE(s)=\EE(s_0)+V,
$$
where $s_0:=s|_{U\times U}$ (note that $\EE(s_0)=\EE(s)\cap U$).
Let $\Vol(s_0)$ denote the volume of $\EE(s_0)$ with respect
to the Lebesgue measure of $U$. Since $I_0$ is bounded,
there exists an $s\in\CF$ with $\Vol(s_0)<+\infty$. Observe that for any
basis $e=(e_1,\dots,e_m)$  of $U$  ($m:=\dim_{\CC}U$) we have
$$
\Vol(s_0)=\frac{C(e)}{\det S},
$$
where $C(e)>0$ is a constant (independent of $s$) and $S=S(s_0)$ denotes
the matrix representation of $s_0$ in the basis $e$,
i.e. $S_{j,k}:=s(e_j,e_k)$, $j,k=1,\dots,m$. In particular, if $U=\CC^m\times
\{0\}^{n-m}$ and $e=(e_1,\dots,e_m)$ is the canonical basis, then
$C(e)$ is the volume of the open unit Euclidean ball $\BB_m\subset\CC^m$.
 We are interested in finding an $s\in\CF$ for
which $\Vol(s_0)$ is minimal, equivalently, $\det S(s_0)$ is maximal.
Observe that if $s$ has this property (with respect to $h$), then for
any $\CC$--linear isomorphism $L:\CC^n\too\CC^n$, the scalar product
$$
\CC^n\times\CC^n\ni(X,Y)\overset{L(s)}\too s(L(X),L(Y))\in\CC
$$
has the extremal property with respect to $h\circ L$. In particular, this permits
us to reduce the situation to the case where $U=\CC^m\times
\{0\}^{n-m}$ and next (by restricting all the above objects to
$\CC^m\simeq\CC^m\times\{0\}^{n-m}$) to assume that $m=n$.

\begin{theorem}[\cite{Wu1987}, \cite{Wu1993}]\label{exactlyone}
There exists exactly one element $s^h\in\CF$ such that
$$
\Vol(s^h_0)=\min\{\Vol(s_0): s\in\CF\}<+\infty.
$$
\end{theorem}

Put $\wdht s^h:=m\cdot s^h$ ($m:=\dim U(h)$), $\WW h:=q_{\wdht s^h}$ .
Obviously, $\WW h\leq\sqrt{m}h$ and $\WW h\equiv\sqrt{m}h$
iff $h=q_s$ for some pseudo--Hermitian scalar product $s$. For instance,
$\WW\|\;\|=\sqrt{n}\|\;\|$, where $\|\;\|$ is the Euclidean
norm in $\CC^n$. Moreover, $\WW(\WW h)\equiv\sqrt{m}\WW h$.

\begin{remark}
Assume that $U(h)=\CC^n$. Let $L:\CC^n\too\CC^n$ be a $\CC$--linear isomorphism such that
$|\det L|=1$ and $h\circ L=h$. Then $\Vol(s^h)=\Vol(L(s^h))$
and hence $s^h=L(s^h)$, i.e.~$s^h(X,Y)=s^h(L(X),L(Y))$, $X,Y\in\CC^n$.
\end{remark}

\begin{theorem}[\cite{Wu1987}, \cite{Wu1993}\footnote{See also
M.~Jarnicki,  P.~Pflug, \textit{Invariant distances and metrics in complex
analysis --- revisited}, to appear.}]\label{pointwise}
{\rm (a)} $h\leq\WW h\leq\sqrt{m}h$.

{\rm (b)} If $h(X):=\max\{h_1(X_1),\;h_2(X_2)\}$, $X=(X_1,X_2)\in
\CC^{n_1}\times\CC^{n_2}$, then
$$
\wdht s^h(X,Y)=\wdht s^{h_1}(X_1,Y_1)+\wdht s^{h_2}(X_2,Y_2),\quad
X=(X_1,X_2),\;Y=(Y_1,Y_2)\in\CC^{n_1}\times\CC^{n_2},
$$
In particular,
$$
\WW h(X)=\Big((\WW h_1(X_1))^2+(\WW h_2(X_2))^2\Big)^{\!1/2},
\quad X=(X_1,X_2).
$$
\end{theorem}

For a domain $G\subset\CC^n$, let $\EM(G)$
denote the space of all {\it pseudometrics}
$$
\eta: G\times\CC^n\too\RR_+,\quad \eta(a;tX)=|t|\eta(a;X),\quad
a\in G,\;X\in\CC^n,\;t\in\CC,
$$
such that
$$
\forall_{a\in G}\;\; \exists_{M, r>0}:\; \eta(z;X)\leq
M\|X\|,\quad z\in\BB(a,r)\subset G,\; X\in\CC^n,
$$
where $\BB(a,r)$ is the open Euclidean ball centered at $a$
with radius $r$.

For $\eta\in\EM(G)$ we define the {\it Wu pseudometric}
$$
(\WW\eta)(a;X):=(\WW\wdht\eta(a;\cdot))(X),\quad a\in
G,\;X\in\CC^n,
$$
where
$$
\wdht\eta(a;X):=\sup\{h(X): \; h \text{ is a $\CC$--seminorm},\;
h\leq\eta(a;\cdot)\},\quad a\in G,\; X\in\CC^n,
$$
denotes the {\it Buseman pseudometric} associated to $\eta$
(cf.~\cite{JarPfl1993b}, \S\,4.3). Observe that $\WW\eta\in\EM(G)$.

Recall that an upper semicontinuous metric $\eta\in\EM(G)$ is said to be
{\it complete} if any $\int\!\eta$--Cauchy sequence is convergent to a point
from $G$, where $\int\!\eta$ denotes the {\it integrated form of $\eta$}
(cf.~\cite{JarPfl1993b}, \S\S\;4.3, 7.3).

\begin{proposition}\label{PropOfWu}
{\rm (a)} If $\eta\in\EM(G)$ is a continuous metric, then
so is $\WW\eta$ (cf.~Example \ref{ex1}).

{\rm (b)} If $\eta\in\EM(G)$ is a continuous complete metric,
then so is $\WW\eta$.

{\rm (c)} If $(\delta_G)_G$ is a holomorphically contractible
family of pseudometrics (cf. \linebreak \cite{JarPfl1993b}),
then for any biholomorphic mapping $F:G\too D$ ($G, D\subset\CC^n$) we have
$$
(\WW\delta_D)(F(z);F'(z)X)=(\WW\delta_G)(z;X),\quad z\in G,\; X\in\CC^n.
$$

{\rm (d)} If $(\delta_G)_G$ is a holomorphically contractible
family of pseudometrics, then for any holomorphic mapping $F:G\too D$
($G\subset\CC^{n_1}$, $D\subset\CC^{n_2}$) we have
$$
(\WW\delta_D)(F(z);F'(z)X)\leq\sqrt{n_2}(\WW\delta_G)(z;X),\quad z\in G,\;
X\in\CC^n,
$$
but, for example, the family $(\WW\kappa_G)_G$ is not holomorphically
contractible, where $\kappa_G$ is the {\it Kobayashi--Royden
pseudometric of $G$} (cf.~Example \ref{ex0}).
\end{proposition}

In the case $\eta=\kappa_G$, the above properties (a) --- (d) were formulated
(without proof) in \cite{Wu1987}, \cite{Wu1993}.

\begin{proof} (a) Fix a point $z_0\in G\subset\CC^n$.
Let $s^z:=s^{\eta(z;\cdot)}$, $z\in G$. We are going to show that
$s^z\underset{z\to z_0}\too s^{z_0}$.

By our assumptions, there exist $r>0$, $c>0$ such that
$$
\eta(z;X)\geq c\|X\|,\quad z\in\BB(z_0,r)\subset G,\; X\in\CC^n.
$$
In particular, the sets
$$
I(z):=\{X\in\CC^n:\wdht\eta(z;X)<1\},\quad z\in\BB(z_0,r),
$$
are contained in the ball $\BB(0,C)$ with $C:=1/c$. Moreover,
$$
|\eta(z;X)-\eta(z_0;X)|\leq\phi(z)\|X\|,\quad X\in\CC^n,
$$
where $\phi(z)\underset{z\to z_0}\too0$. Hence
$$
(1+C\phi(z))^{-1} I(z)\subset I(z_0)\subset(1+C\phi(z))I(z),\quad z\in\BB(z_0,r),
$$
and consequently,
\begin{align*}
I(z_0)\subset(1+C\phi(z))\EE(s^z)&=\EE((1+C\phi(z))^{-2}s^z), \tag{*}\\
I(z)\subset(1+C\phi(z))\EE(s^{z_0})&=\EE((1+C\phi(z))^{-2}s^{z_0}),\quad
z\in\BB(z_0,r).
\end{align*}
Hence
\begin{align*}
\Vol(s^{z_0})\leq\Vol((1+C\phi(z))^{-2}s^z)&=(1+C\phi(z))^{2n}\Vol(s^z),\\
\Vol(s^z)\leq\Vol((1+C\phi(z))^{-2}s^{z_0})&=(1+C\phi(z))^{2n}\Vol(s^{z_0}),\quad
z\in\BB(z_0,r).
\end{align*}
Thus $\Vol(s^z)\underset{z\to z_0}\too\Vol(s^{z_0})$.

Take a sequence $z_\nu\too z_0$. Since
$$
|s^{z_\nu}(e_j,e_k)|\leq\eta(z_\nu;e_j)\eta(z_\nu;e_k),\quad
j,k=1,\dots,n,\;\nu\in\NN,
$$
we may assume that $s^{z_\nu}\too s^\ast$, where $s^\ast$ is a pseudo--Hermitian scalar
product. We already know that $\Vol(s^\ast)=\Vol(s^{z_0})$. Moreover, by
(*), $I(z_0)\subset\EE(s^\ast)$. Consequently, the uniqueness of $s^{z_0}$ implies
that $s^\ast=s^{z_0}$.

(b) Recall that $\int\eta=\int\wdht\eta$ --- cf.~\cite{JarPfl1993b},
Proposition 4.3.5(b). By (a),
$\WW\eta$ is a continuous metric. In particular, the
distance $\int\!(\WW\eta)$ is well defined.
By Theorem \ref{pointwise}(a) we get
$$
\int\!\wdht\eta\leq\int\!(\WW\eta),
$$
which directly implies the required result.

(c) The result is obvious because for any $z\in G$, the mapping
$F'(z)$ is a $\CC$--linear isomorphism and
$\delta_D(F(z);F'(z)X)=\delta_G(z;X)$, $X\in\CC^n$.

(d) It is known that the family $(\wdht\delta_G)_G$ is holomorphically
contractible (\cite{JarPfl1993b}, Theorem 4.3.10(c)). Hence, using
Theorem
\ref{pointwise}(a), we get
\begin{multline*}
(\WW\delta_D)(F(z);F'(z)X)\leq\sqrt{n_2}\wdht\delta_D(F(z);F'(z)X)\\
\leq\sqrt{n_2}\wdht\delta_G(z;X)\leq\sqrt{n_2}(\WW\delta_G)(z;X),\quad
 z\in G,\;X\in\CC^n.
\end{multline*}
\end{proof}

\begin{example}\label{ex0}
Let $G:=\{(z_1,z_2)\in\BB_2: |z_1|<\eps\}$, $0<\eps<1/\sqrt{2}$.
Recall that $\kappa_{\BB_2}(0;X)=\|X\|$ and
$\kappa_G(0;X)=\max\{\|X\|, |X_1|/\eps\}$, $X=(X_1,X_2)$.
Then $(\WW\kappa_{\BB_2})(0;(0,1))=\sqrt{2}>(\WW\kappa_G)(0;(0,1))$.
In particular, the family $(\WW\kappa_D)_D$ is not contractible
with respect to inclusions.
\end{example}

We point out that Proposition \ref{PropOfWu}(a) gives us the
continuity of $\WW\eta$ only in the case where $\eta$ is a
continuous metric. It is natural to conjecture that in the general
case, where $\eta$ is only an upper semicontinuous (pseudo)metric,
$\WW\eta$ remains to be upper semicontinuous.
The following Example \ref{ex1} shows that in general this is not
true.
In the case where $\eta$ is a continuous pseudometric, we do not
know whether $\WW\eta$ is upper semicontinuous. Observe that
the upper semicontinuity (or at least Borel measurability) of $\WW\eta$
appears in a natural way when one defines $\int\!(\WW\eta)$.
In the case where
$\eta=\kappa_G$, the upper semicontinuity of $\WW\kappa_G$ is
claimed for instance in \cite{Wu1993} (Theorem 1), \cite{CheKim1996}
(Proposition 2), \cite{Juc2003} (Theorem 0), but so far there is no proof.

\begin{example}\label{ex1}
There is  an upper semicontinuous metric $\eta$
such that $\WW\eta$ is not upper semicontinuous.

Indeed, let $\eta:\BB_2\times\CC^2\too\RR_+$, $\eta(z;X):=\|X\|$ for $z\neq0$,
and $\eta(0,X):=\max\{\|X\|,|X_1|/\eps\}$, $X=(X_1,X_2)\in\CC^2$ ($\eps>0$ small).
Then $(\WW\eta)(z;X)=\sqrt{2}\|X\|$ for $z\neq0$, and
$\{X\in\CC^2: (\WW\eta)(0;X)<1\}\setminus\BB(0,1/\sqrt{2})\neq\varnothing$, so $\WW\eta$
is not upper semicontinuous (cf.~Example \ref{ex0}).
\end{example}

\begin{example}\label{ex3}
There exists a bounded domain $G\subset\CC^2$ such that $\WW\kappa_G$
is not continuous (see Proposition 2 in \cite{CheKim1996}, where such a
continuity is claimed).

Indeed,
let $D\subset\CC^2$ be a domain such that (cf.~\cite{JarPfl1993b},
Example 3.5.10):

$\bullet$ there exists a dense subset $M\subset\CC$ such that $(M\times\CC)\cup
(\CC\times\{0\})\subset D$,

$\bullet$ $\kappa_D(z;(0,1))=0$, $z\in A:=M\times\CC$,

$\bullet$ there exists a point $z^0\in D\setminus A$ such that $\kappa_D(z^0;X)
\geq c\|X\|$, $X\in\CC^2$, where $c>0$ is a constant.

For $R>0$ let $D_R:=\{z=(z_1,z_2)\in D: |z_j-z^0_j|<R,\;j=1,2\}$. It is known
that $\kappa_{D_R}\searrow\kappa_D$ when $R\nearrow+\infty$. Observe that
$z^0\in D_R$ and
$$
\kappa_{D_R}(z^0;X)\geq\kappa_D(z^0;X)\geq c\|X\|,\quad X\in\CC^2.
$$
Hence, by Theorem \ref{pointwise}(a), $(\WW\kappa_{D_R})(z^0;X)\geq
c\|X\|$, $X\in\CC^2$. In particular,
$$
(\WW\kappa_{D_R})(z^0;(0,1))\geq c.
$$
Fix a sequence $M\ni z_k\too z^0_1$.
Note that $\{z_k\}\times (z^0_2+RE)\subset D_R$, which implies that
$\kappa_{D_R}((z_k,z^0_2);(0,1))\leq1/R$, $k=1,2,\dots$. In particular,
$$
(\WW\kappa_{D_R})((z_k,z^0_2);(0,1))\leq
\sqrt{2}\kappa_{D_R}((z_k,z^0_2);(0,1))\leq\sqrt{2}/R,\quad k=1,2,\dots.
$$

Now it clear that if $R>\frac{\sqrt{2}}{c}$, then
$$
\limsup_{k\to+\infty}
(\WW\kappa_{D_R})((z_k,z^0_2);(0,1))\leq\sqrt{2}/R<c\leq
(\WW\kappa_{D_R})(z^0;(0,1)),
$$
which shows that for $G:=D_R$ the pseudometric $\WW\kappa_G$ is not continuous.
\end{example}

\begin{remark}
We point out the influence of the factor $\sqrt{m}$ in the definition of
$\WW$ to its upper semicontinuity.

Suppose we defined $\wdtl\WW h:=q_{s^h}$,
$\wdtl\WW\eta(a;X)=(\wdtl\WW\wdht\eta(a;\cdot))(X)$.
Then, using the product formula (Theorem \ref{pointwise}(b)), we would get
a domain $G\subset\CC^3$ such
that $\wdtl\WW\kappa_G$ is not upper semicontinuous.

Indeed (the example is due to W.~Jarnicki), let $D\subset\CC^2$
and $D\ni z_k\too z_0\in D$ be such that:

$\bullet$ $\kappa_D(z_k;\cdot)$ is not a metric (in particular,
$m(k):=\dim U(\wdht\kappa_D(z_k;\cdot))\leq 1$, $k\in\NN$),

$\bullet$ $\kappa_D(z_0;\cdot)$ is a metric.

Take, for instance, the domain $D$ as in Example \ref{ex3}.

Put $G:=D\times E\subset\CC^3$.
Then $\wdtl\WW\kappa_G$ is not upper semicontinuous at $((z_0,0),(0,1))$
because
\begin{align*}
(\wdtl\WW\kappa_G)^2((z_k,0);(0,1))&=s^{\kappa_G((z_k,0);\cdot)}((0,1),(0,1))
=\frac1{m(z_k)+1}\geq\frac12,\quad k\in\NN,\\
(\wdtl\WW\kappa_G)^2((z_0,0);(0,1))&=s^{\kappa_G((z_0,0);\cdot)}((0,1),(0,1))
=\frac1{m(z_0)+1}=\frac13.
\end{align*}
\end{remark}

We conclude the paper by repeating the main open question.

{\bf Problem.} Let $\eta\in\{\gamma_G^{(k)}, A_G, \kappa_G\}$ (cf.~\cite{JarPfl1993b}).
Is $\WW\eta$ upper semicontinuous~?

\bibliographystyle{amsplain}

\begin{thebibliography}{XXXXXXXXX}
\makeatletter\renewcommand{\@biblabel}[1]{[#1]}\makeatother
\bibitem[Che-Kim 1996]{CheKim1996}
C.~K.~Cheung,  K.~T.~Kim,
\textit{Analysis of the Wu metric. I: The case of convex Thullen domains},
Trans. Amer. Math. Soc. 348 (1996), 1421--1457.
\bibitem[Che-Kim 1997]{CheKim1997}
C.~K.~Cheung,  K.~T.~Kim,
\textit{Analysis of the Wu metric. II: The case of non-convex Thullen domains},
Proc. Amer. Math. Soc. 125 (1997), 1131--1142.
\bibitem[Che-Kim 2003]{CheKim2003}
C.~K.~Cheung,  K.~T.~Kim,
\textit{The constant curvature property of the Wu invariant metric},
Pac. J. Math. 211 (2003), 61--68.
\bibitem[J-P 1993]{JarPfl1993b}
M.~Jarnicki,  P.~Pflug,
\textit{Invariant Distances and Metrics in Complex Analysis},
 de Gruyter Expositions in Mathematics 9, Walter de Gruyter  1993.
\bibitem[Juc 2002]{Juc2003}
P.~Jucha,
\textit{The Wu metric in elementary Reinhardt domains},
Univ. Iag. Acta Math. 40 (2002), 83--89.
\bibitem[Kim 1998]{Kim1998}
K.~T.~Kim,
\textit{The Wu metric and minimum ellipsoids}, Proc. 3rd Pacific Rim
Geometry Conf. (J.~Choe ed.),
Monogr. Geom. Topology 25, International Press, Cambridge, 1998, 121--138.
\bibitem[Wu]{Wu1987}
H.~Wu,
\textit{Unpublished notes}.
\bibitem[Wu 1993]{Wu1993}
H.~Wu,
\textit{Old and new invariant metrics},
Several complex variables: Proc. Mittag--Leffler Inst. 1987--88
(J.E.~Forn{\ae}ss ed.), Math. Notes, Princeton Univ. Press 38 (1993), 640--682.
\end{thebibliography}

\end{document}